\documentclass{article}
\PassOptionsToPackage{numbers,sort&compress}{natbib}

\usepackage[final]{neurips_2025_ml4ps}

\makeatletter
\def\NAT@warning#1{}
\makeatother

\usepackage[utf8]{inputenc} % allow utf-8 input
\usepackage[T1]{fontenc}    % use 8-bit T1 fonts
\usepackage{hyperref}       % hyperlinks
\usepackage{url}            % simple URL typesetting
\usepackage{booktabs}       % professional-quality tables
\usepackage{amsfonts}       % blackboard math symbols
\usepackage{nicefrac}       % compact symbols for 1/2, etc.
\usepackage{microtype}      % microtypography
\usepackage{xcolor}         % colors
\usepackage{subcaption}
\usepackage{graphicx}
\usepackage{amsmath}

\title{Offline Maximizing Minimally Invasive Proper Orthogonal Decomposition for Reduced Order Modeling of $S_n$ Radiation Transport}

\author{%
  Quincy Huhn\\
  Department of Nuclear Engineering\\
  Texas A\&M University\\
  College Station, TX 77584\\
  \texttt{quincy.huhn98@tamu.edu} \\
  \And
  Jean C. Ragusa\\
  Department of Nuclear Engineering\\
  Texas A\&M University\\
  College Station, TX 77584\\
  \texttt{jean.ragusa@tamu.edu} \\
  \And
  Youngsoo Choi\\
  Center for Applied Scientific Computing\\
  Lawrence Livermore National Laboratory\\
  Livermore, CA 94550\\
  \texttt{choi15@llnl.gov}
}

\begin{document}

\maketitle

\begin{abstract}
    Deterministic solutions to the $S_n$ transport equation can be computationally expensive to calculate. Reduced Order Models (ROMs) provide an efficient means of approximating the Full Order Model (FOM) solution. We propose a novel approach for constructing ROMs of the $S_n$ radiation transport equation, Offline Maximizing Minimally Invasive (OMMI) Proper Orthogonal Decomposition (POD). POD uses snapshot data to build a reduced basis, which is then used to project the FOM. Minimally Invasive POD leverages the sweep infrastructure within deterministic $S_n$ transport solvers to construct the reduced linear system, even though the FOM linear system is never directly assembled.
    OMMI-POD extends Minimally Invasive POD by performing transport sweeps offline, thereby maximizing the potential speedup. It achieves this by generating a library of reduced systems from a training set, which is then interpolated in the online stage to provide a rapid approximate solution to the $S_n$ transport equation.
    The model’s performance is evaluated on a multigroup 2-D test problem, demonstrating low error and a 1600-fold speedup over the full order model. 
\end{abstract}

\section{Introduction}
The steady-state transport equation is six-dimensional (three in space, two in angle, and one in energy), so obtaining solutions for high-fidelity problems can become very expensive. This is particularly true in multi-query settings such as design optimization and uncertainty quantification, where many similar problems must be solved repeatedly. Reduced Order Models (ROMs) can be helpful in these cases, as they provide approximate solutions to the Full Order Model (FOM) in significantly less time.

The ROMs developed in this work are constructed using Proper Orthogonal Decomposition (POD)\cite{Pearson01111901}, which has been widely applied to problems in nuclear engineering \cite{GERMAN2019144,HUGHES,TENCER}. However, typical POD methods are challenging to apply to transport problems due to the size of the system. For finely discretized problems, the assembly of the reduced system could either be prohibitively expensive to perform online or infeasible, depending on the available computing resources. Current approaches to applying POD to transport problems either address this issue or leverage the POD basis for other aspects of the problem. For example, Space-time ROM utilizes POD to reduce both time and space in transient transport problems \cite{CHOI2021109845}. Affine decomposition has been used to shift the expensive multiplication of the system matrix by the modes offline, enabling efficient online assembly of the reduced system \cite{doi:10.1080/00295639.2022.2112901}. POD has also been used to construct optimal angular basis functions for transport problems \cite{BUCHAN2015138}, accelerate Monte Carlo simulations \cite{UDAGEDARA2015237}, approximate the Eddington tensor in moment calculations \cite{COALE2024113044}, and form ROMs for a finite volume formulation of transport \cite{SUN2020107799}.

There are two types of ROMs: intrusive and non-intrusive. Non-intrusive ROMs are purely data-driven, relying either on data generated from solutions to the underlying PDE or from physical observations. By not requiring access to the PDE, non-intrusive ROMs can be readily applied to a wide range of problems. Intrusive ROMs, on the other hand, frequently require significant modifications to the FOM to implement. As a result, intrusive ROMs are often more accurate but typically require extensive initial setup for a new application. 

Minimally Invasive POD \cite{BEHNE2022111525} is an intrusive ROM and therefore benefits from increased accuracy, but it is also designed to require minimal modification to the FOM code. This is possible due to the typical solution method that is used when solving radiation transport problems. Instead of assembling the full system matrix, transport codes will use a matrix-free method known as a sweep to perform the operator action of the system matrix. The operator action of the system matrix is supplied to the solver instead of the system matrix itself, and Minimally Invasive POD uses this operator action to assemble the ROM. Since the operator action is typically straightforward to expose to external programs, Minimally Invasive POD can function as an intrusive ROM while still requiring minimal modifications to typical transport codes.

This paper extends Minimally Invasive POD to create Offline Maximizing Minimally Invasive POD (OMMI-POD), a novel approach for building a ROM for the transport equation. OMMI-POD maximizes the offline phase by performing all required sweeps offline and using interpolation during the online stage.

\section{Methods}
\subsection{Neutron Transport FOM}
The steady-state neutron transport equation is
\begin{align}\label{eq:BTE} 
        \vec{\Omega} \cdot \vec{\nabla}\psi(\vec{r},\vec{\Omega},E)  + \sigma_t(\vec{r},E) \psi(\vec{r},\vec{\Omega},E)
     =& \int_{4\pi} d\vec{\Omega}'\int_0^{\infty}dE' \sigma_s(E' \rightarrow E, \vec{\Omega}' \rightarrow \vec{\Omega}) \psi(\vec{r},\vec{\Omega}',E')  \nonumber \\&+ Q_{e}(\vec{r},\vec{\Omega},E) 
\end{align}
where $\psi$ is the angular flux, $\vec{r}$ is the position, $E$ is the particle energy, $\vec{\Omega}$ is the direction. The total cross section is denoted by $\sigma_t$, and $\sigma_s(E' \rightarrow E, \vec{\Omega}' \rightarrow \vec{\Omega})$ is the differential scattering cross section from direction $\vec{\Omega}'$ to $\vec{\Omega}$ and from energy $E'$ to $E$. In this paper, the transport equation is discretized in angle using the discrete ordinates ($S_N$) method, in energy using the multigroup method, and in space using Discontinuous Galerkin finite elements.

The fully discretized transport equation can be written in operator form as
\begin{equation}
    \mathsf{L}\psi = \mathsf{M}\mathsf{S}\phi + \mathsf{Q}_e,
\end{equation}
where $\mathsf{L}$ is the loss operator which combines the streaming and interaction terms, $\mathsf{M}$ is the moment to discrete operator, and $\mathsf{S}$ is the scattering operator. The scalar flux is given by $\phi = \mathsf{D}\psi$, where $\mathsf{D}$ is the discrete to moment operator, and $\mathsf{Q}_e$ is the discretized external source. 

To solve the transport equation, it is reformulated as the system
\begin{equation}\label{eq:GMRES}
    (\mathsf{I} - \mathsf{D}\mathsf{L}^{-1}\mathsf{M}\mathsf{S})\phi = \mathsf{D}\mathsf{L}^{-1}\mathsf{Q}_e
\end{equation}
which is solved using GMRES. GMRES is used for transport problems because of the large size of typical transport systems. The six-dimensional nature of the $S_N$ transport equation leads to systems so large that they cannot be stored explicitly. GMRES addresses this limitation by not requiring the system matrix to be assembled; only the operator action of the system matrix needs to be performed.

In transport problems, this operator action is performed via the matrix-free sweep. A sweep is a systematic method that inverts the linear system cell by cell, propagating results downstream. By starting at a boundary and advancing through the medium, a sweep effectively recreates the global inversion of $\mathsf{L}$ without requiring the assembly of the system matrix.

\subsection{OMMI-POD ROM}
A projection-based reduced order model for a linear system $Ax=b$ is given by

\begin{equation}
 W^TAUc=W^Tb
\end{equation}
where $U$ is the reduced basis, $W$ is the reduced test space, and $c$ is the vector of expansion coefficients used to approximate the solution as $x \approx Uc$. The reduced test space can be chosen as $W=U$, corresponding to Galerkin projection, or as $W=AU$, corresponding to Petrov-Galerkin projection. For radiation transport problems, Petrov-Galerkin projection is preferred because it guarantees stability of the ROM for hyperbolic systems \cite{BUI-THANH}. 

To obtain the reduced basis, a snapshot matrix $X = [x_0, x_1,\dots,x_{N_{s}}]$ is constructed from the training data. The reduced basis is then obtained by performing the Singular Value Decomposition (SVD) of the snapshot matrix $X = U \Lambda V^T$, where $U$ contains the POD modes and $\Lambda$ is a diagonal matrix of the singular values. The modes in $U$ are truncated to retain a specified amount of information, defined as
\begin{equation}
    I(r) = \frac{\sum_{i=1}^r \lambda_{i}^2}{\sum_{i=1}^{N_s} \lambda_{i}^2}
\end{equation}
where $I(r)$ is the information retained when truncating at rank $r$. The rank $r$ is typically selected to minimize the number of modes retained while ensuring $I(r)$ is greater than some threshold value.

This formulation of the ROM is problematic for radiation transport problems due to the prohibitive size of $A$. As described in the previous section, the matrix is never assembled in transport solvers; instead, a matrix-free sweep is performed. However, the sweep can be used to apply the operator action of $A$ to each mode in the reduced basis, thereby forming $AU$. This approach is called Minimally Invasive POD because the operator action of $A$ is typically already available to solvers, so minimal code modifications are required to apply $A$ to $U$. 

The main limitation of Minimally Invasive POD is that the computational cost of the sweeps must be incurred online. Because each sweep is on the order of the FOM and $r+1$ sweeps are required to form the reduced system ($r$ to perform $AU$ plus a sweep to form $b$), the potential speedup for Minimally Invasive POD is quite limited. Offline Maximizing Minimally Invasive POD (OMMI-POD) solves this by moving the sweeps offline. During the offline phase, Minimally Invasive POD is applied at each training point to construct a library of reduced systems
\begin{equation}
        \mathfrak{A}_r = [A_{r,0},A_{r,1},\dots, A_{r,N_s}], \mathfrak{b}_r = [b_{r,0},b_{r,1},\dots, b_{r,N_s}]
\end{equation}
where, in the Petrov-Galerkin case, $A_r=(AU)^TAU$ and $b_r = (AU)^Tb$. In the online phase, the reduced systems are interpolated to obtain the reduced system at a new parameter value. The interpolation of $A_r$ is performed with matrix manifold interpolation, and both $A_r$ and $b_r$ are interpolated with radial basis functions. The resulting reduced system can then be solved very rapidly to form the ROM solution. In this way, OMMI-POD greatly reduces costs during the online phase, thereby increasing the potential speedup. The effectiveness of OMMI-POD ultimately depends on the choice of interpolation and sampling strategies, which must be carefully chosen to minimize error across the parameter space.

\section{Results}

We evaluate the performance of OMMI-POD using a two-group version of the 2-D checkerboard problem described in \cite{osti_800993}. The domain contains a background scatterer with $\sigma_a = 0$ and embedded absorbing occlusions with $\sigma_s = 0$, driven by a central isotropic source. Two physical parameters are varied. The first, $\theta_1$, is the total cross section in the absorber and is uniformly distributed over $\theta_1 \in [7.5, 12.5]$, identical for both energy groups. The second, $\theta_2$, is the down-scattering ratio, sampled from $\theta_2 \in [0.5, 1.0]$. The scattering cross sections in the scatterer are therefore $\sigma_{s,0\rightarrow 0}=1-\theta_2$, $\sigma_{s,0\rightarrow 1}=\theta_2$, and $\sigma_{s,1\rightarrow 1}=1.0$. The problem is discretized using 1,024 angular directions and 19,600 spatial degrees of freedom, yielding just over 40 million angular unknowns. OpenSn is used as the FOM solver and to perform the sweeps required for Minimally Invasive POD, while libROM is used to perform singular value decomposition and interpolation.

\begin{table}[!htbp]
  \caption{Speedup and Error of OMMI-POD}
  \label{sample-table}
  \centering
  \begin{tabular}{lll}
    \toprule
    $N_{snap}$     & Speedup     & Error \\
    \midrule
    50     & 1590 & $1.2\times10^{-3}$     \\
    100     & 1190 & $1.2\times10^{-4}$ \\
    150     & 800       & $1.3\times10^{-4}$ \\
    \bottomrule
  \end{tabular}
\end{table}

We generate a training set by sampling a varying number of parameter initializations. A reduced basis with rank $r=5$ is then constructed using POD. The relative $L_2$ error and corresponding speedup for different training set sizes are shown in Table \ref{sample-table}. OMMI-POD achieves very high speedup while maintaining low error across all training set sizes.

\begin{figure}[!htbp]
  \centering
  % First row
  \begin{subfigure}{0.3\textwidth}
    \includegraphics[width=\linewidth]{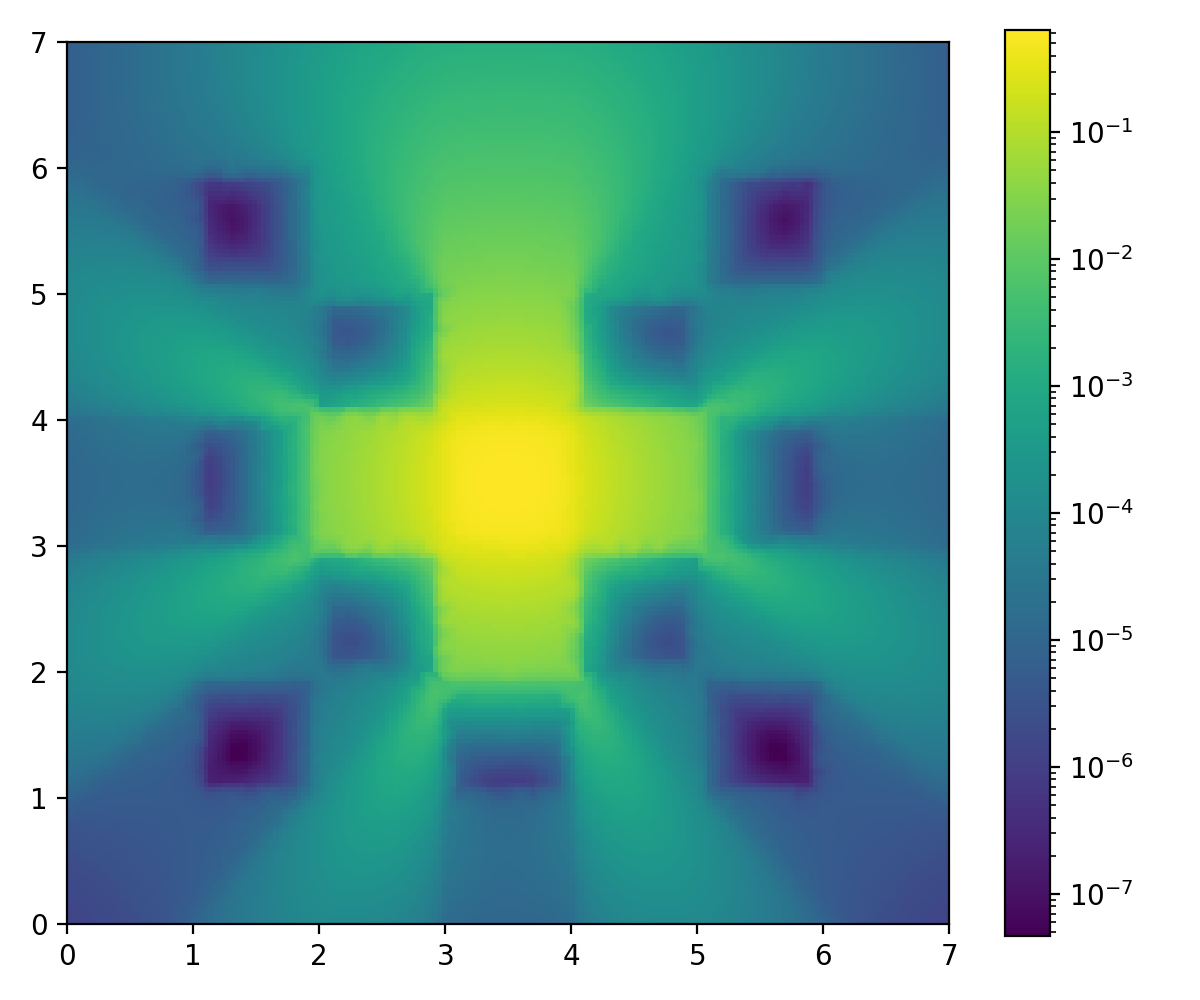}
    \caption{FOM Group 0}
  \end{subfigure}
  \begin{subfigure}{0.3\textwidth}
    \includegraphics[width=\linewidth]{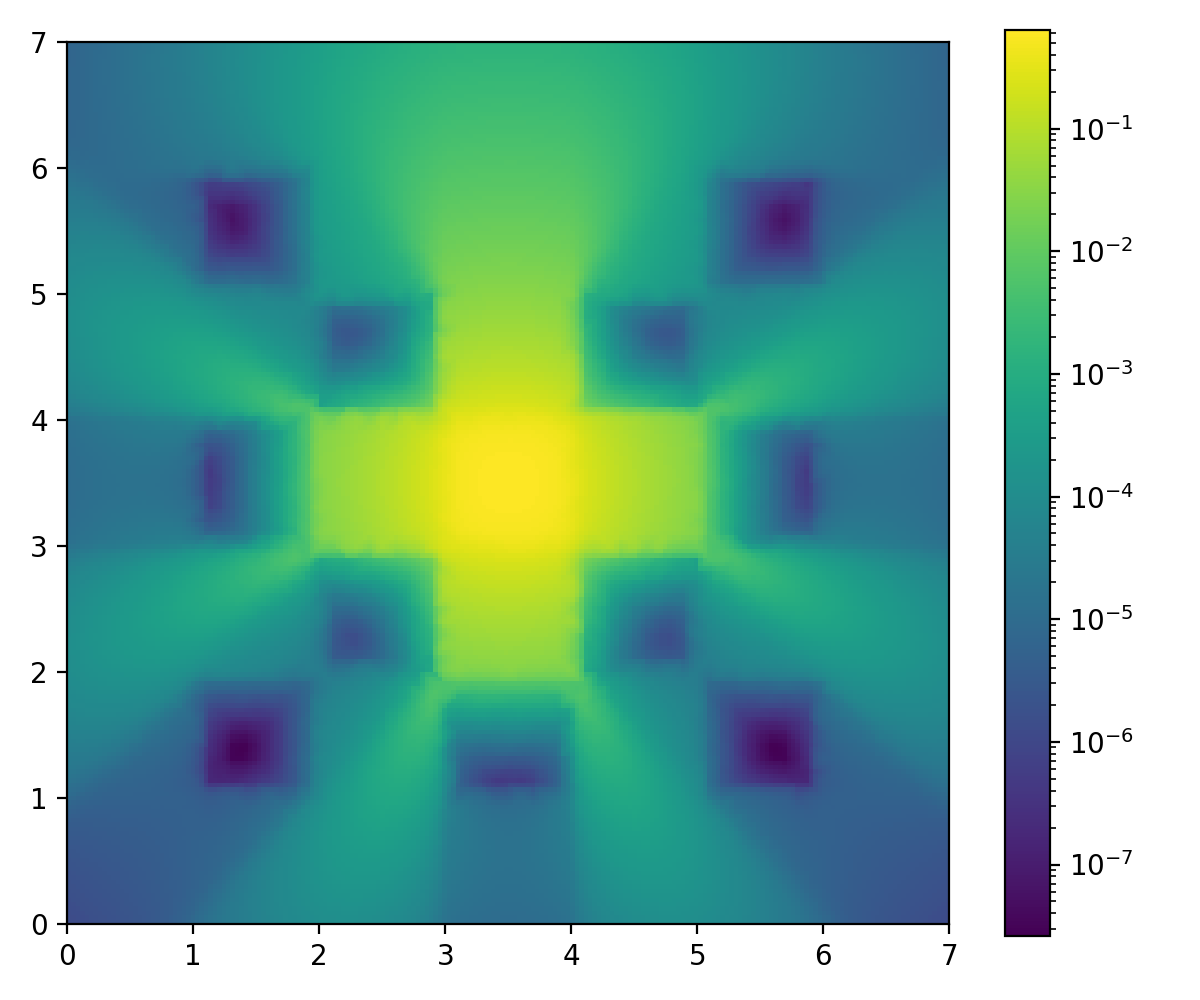}
    \caption{ROM Group 0}
  \end{subfigure}
  \begin{subfigure}{0.3\textwidth}
    \includegraphics[width=\linewidth]{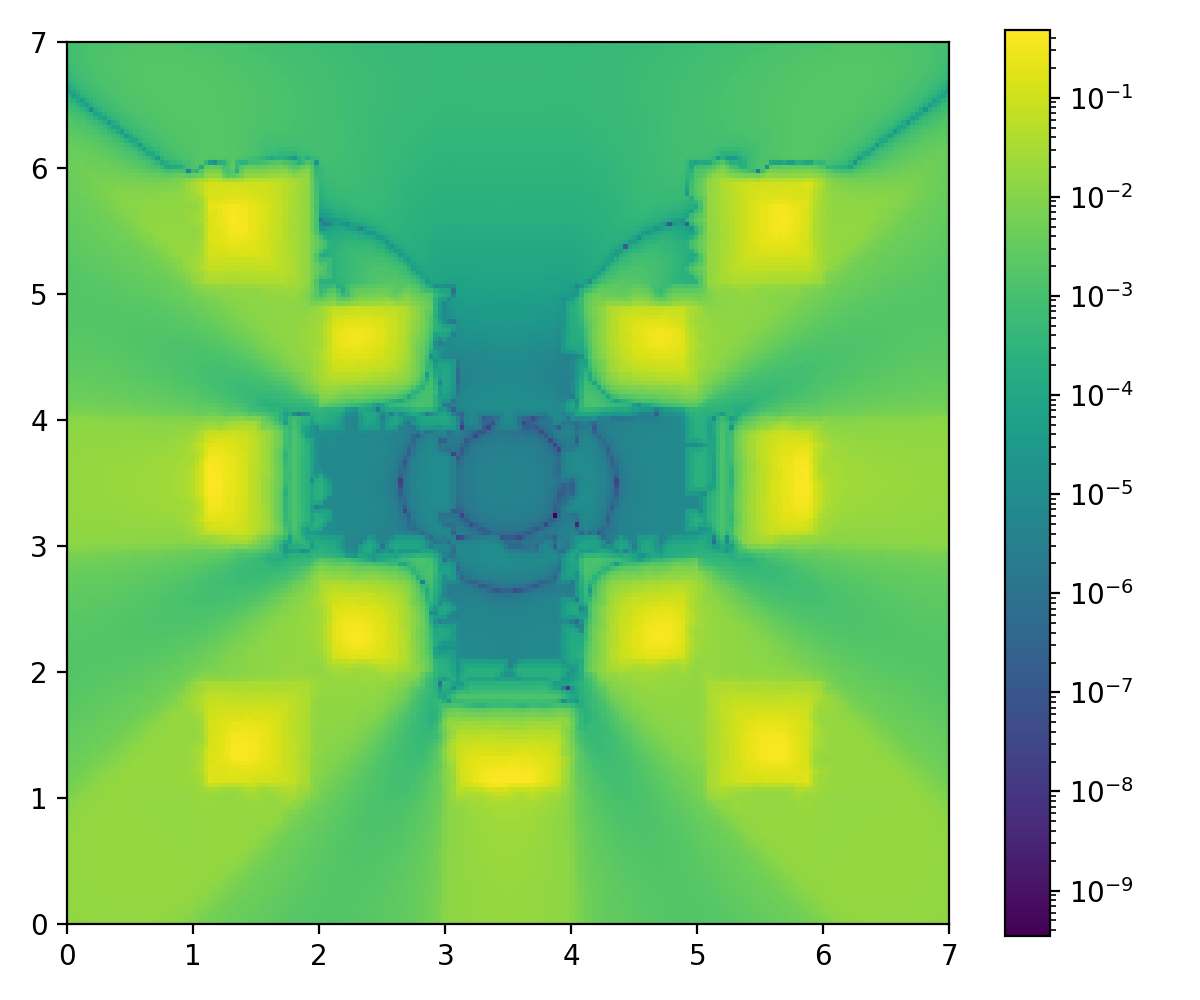}
    \caption{Error Group 0}
  \end{subfigure}

  % Line break between the two sets
  \vspace{0.5cm}

  % Second row
  \begin{subfigure}{0.3\textwidth}
    \includegraphics[width=\linewidth]{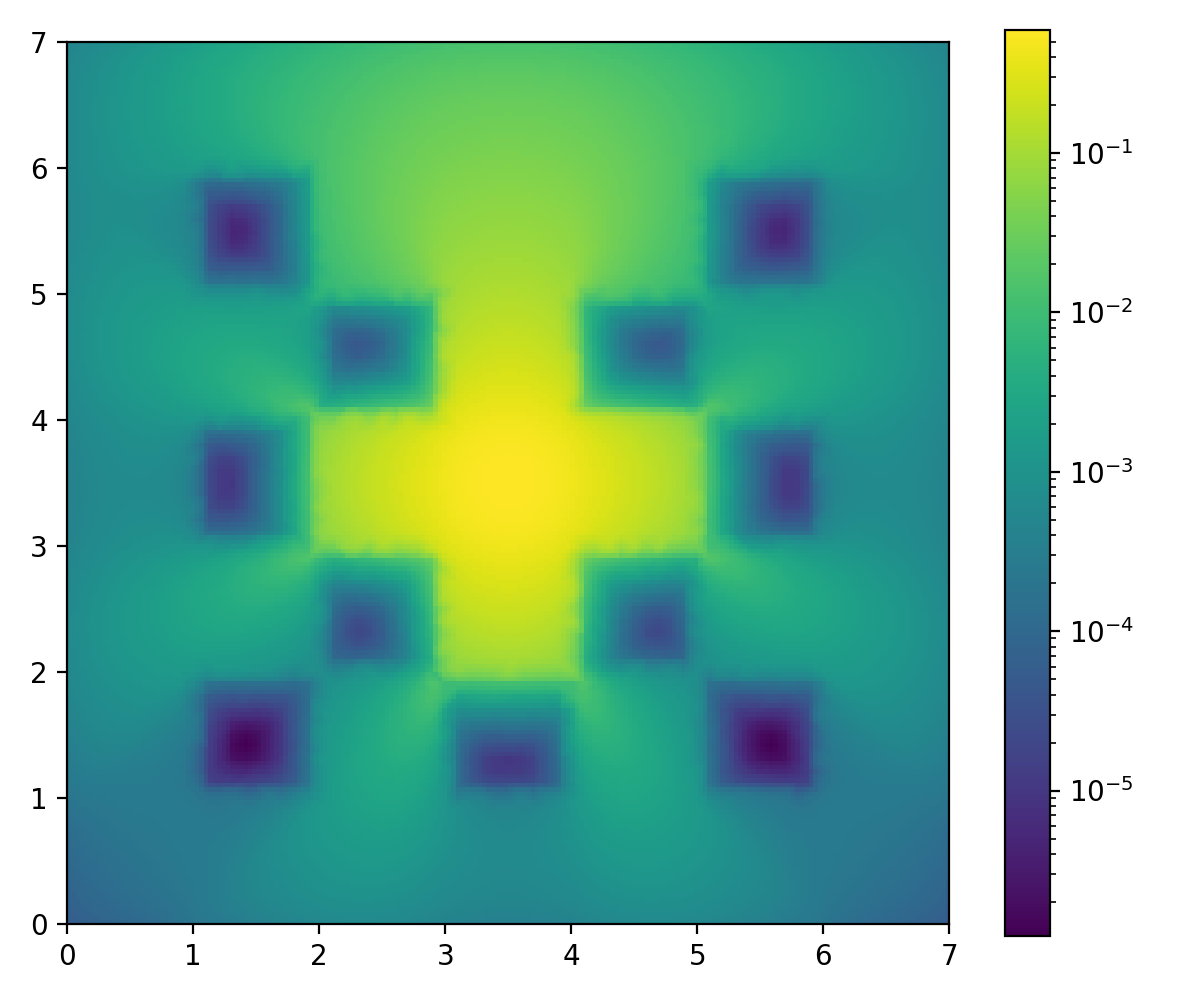}
    \caption{FOM Group 1}
  \end{subfigure}
  \begin{subfigure}{0.3\textwidth}
    \includegraphics[width=\linewidth]{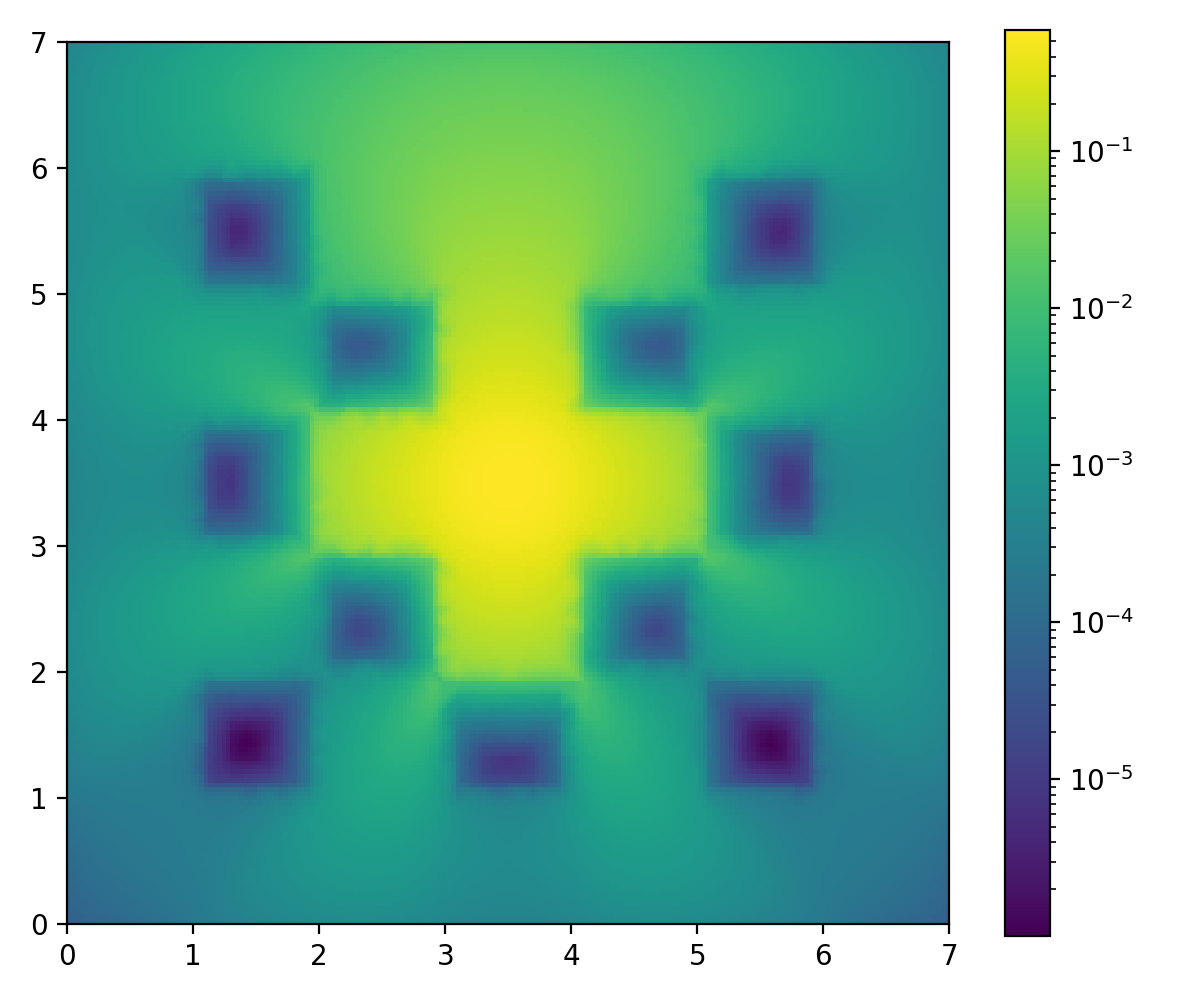}
    \caption{ROM Group 1}
  \end{subfigure}
  \begin{subfigure}{0.3\textwidth}
    \includegraphics[width=\linewidth]{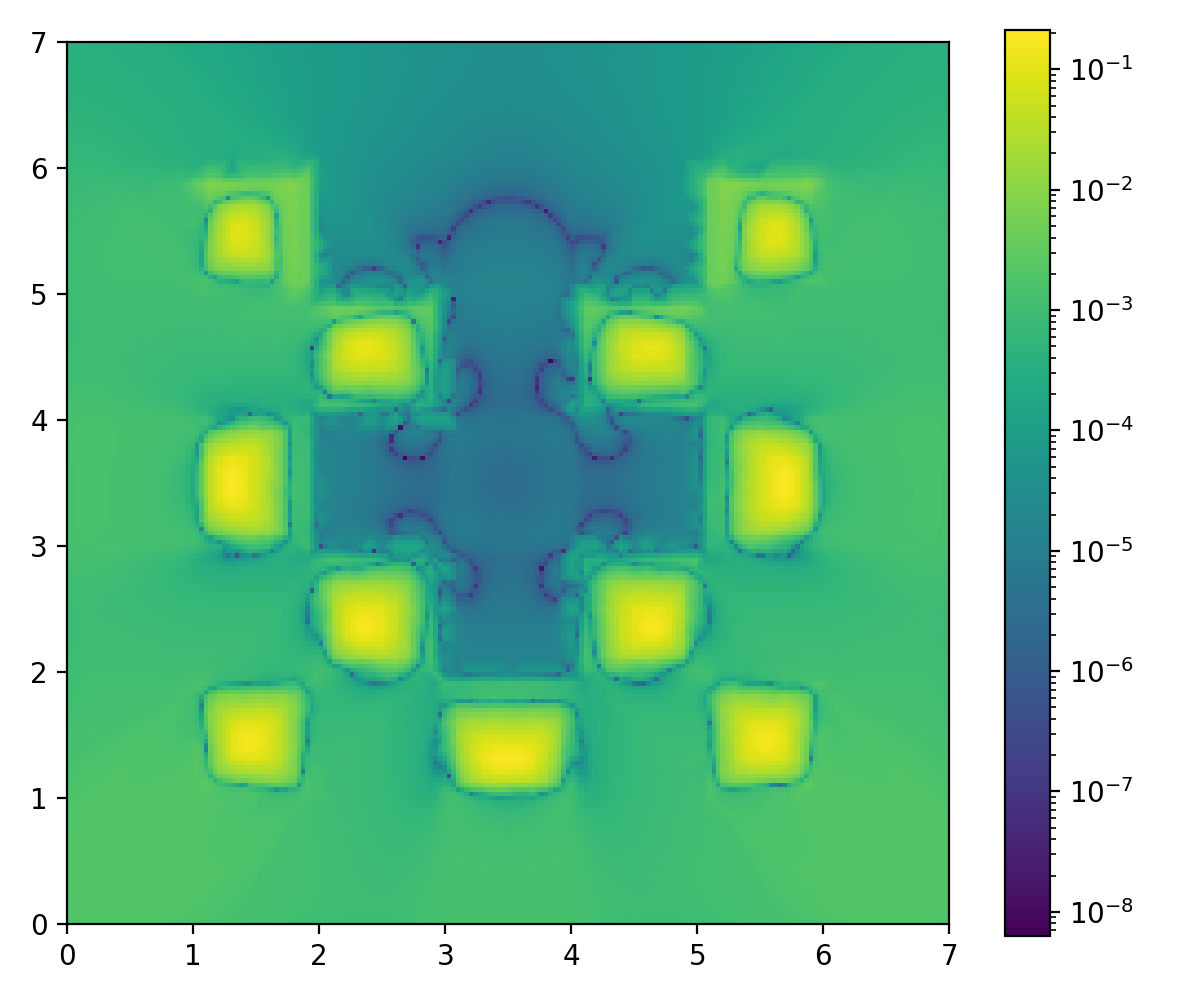}
    \caption{Error Group 1}
  \end{subfigure}

  \caption{ROM and FOM solutions at various parameter initializations.}
  \label{fig:combined}
\end{figure}
Figure \ref{fig:combined} compares the FOM solution, ROM solution, and relative error for both energy groups. The results are generated using 100 snapshots, yielding an average error of $1.2\times10^{-4}$ and an average speedup of 1190 times over a 10-point randomly sampled testing set.

\section{Conclusion}
This paper presents Offline Maximizing Minimally Invasive Proper Orthogonal Decomposition (OMMI-POD), a method of constructing Reduced Order Models (ROMs) for the $S_n$ radiation transport equation. OMMI-POD uses Minimally Invasive POD to create a library of reduced systems, which are then leveraged online to rapidly construct a ROM via interpolation. OMMI-POD achieved an average $L_2$ relative testing error of $1.2\times10^{-4}$ while providing an average speedup of around 1190 times. This represents a substantial improvement over the $1.5$ times speedup of the underlying Minimally Invasive POD. Future research will extend OMMI-POD to $k$-eigenvalue problems and investigate performance improvements through different sampling techniques and interpolation schemes.

\begin{ack}
This work was performed at Texas A\&M under Subcontract No. B661293 with Lawrence Livermore National Security, LLC. Livermore National Laboratory is operated by Lawrence Livermore National Security, LLC, for the U.S. Department of Energy, National Nuclear Security Administration under Contract DE-AC52-07NA27344.
YC was partially supported by the U.S. Department of Energy (DOE), Office of Science, Office of Advanced Scientific Computing Research (ASCR), through the CHaRMNET Mathematical Multifaceted Integrated Capability Center (MMICC) under Award Number DE-SC0023164 and the LEADS SciDAC Institute under Project Number SCW1933. LLNL document release number: LLNL-PROC-2012525.
\end{ack}

% \bibliographystyle{plainnat}
% \bibliography{Bibliography}

%%%%%%%%%%%%%%%%%%%%%%%%%%%%%%%%%%%%%%%%%%%%%%%%%%%%%%%%%%%%

% \appendix

\end{document}